# Optimization in Bochner Spaces

Shuting Ai

Academy of Mathematics and Systems Science
Chinese Academy of Sciences
Zhongguancun, Beijing 100190
China

Jinlu Li

Department of Mathematics
Shawnee State University
Portsmouth Ohio 45662
USA

## Abstract

In this paper, we study some optimization problems in uniformly convex and uniformly smooth Bochner spaces. We consider four cases of the underlying subsets: closed and convex subsets, closed and convex cones, closed subspaces and closed balls. In each case, we study the existence problems and the inverse image properties. We think that the results and the analytic methods in this paper could be applied to the theories of stochastic optimization, stochastic variational inequality and stochastic fixed point.

Keywords: optimization; uniformly convex and uniformly smooth Bochner spaces; normalized duality mapping; metric projection; generalized metric projection

Mathematics subject classification (2010): 65K10, 90C25, 90C26, 90C48

## 1. Introduction

Let $(S, \mathcal{A}, \mu)$ be a positive and complete measure space. Let $(X, \|\cdot\|_X)$ be a real uniformly convex and uniformly smooth Banach space with topological dual space $(X^*, \|\cdot\|_{X^*})$. For $1 < p, q < \infty$ with $\frac{1}{p} + \frac{1}{q} = 1$, let $(L_p(S; X), \|\cdot\|_{L_p(S;X)})$ be the Lebesgue-Bochner function space (simply called the Bochner space). $(L_p(S; X), \|\cdot\|_{L_p(S;X)})$ is a uniformly convex and uniformly smooth Banach space with dual space $(L_p(S; X))^* = L_q(S; X^*)$. $L_p(S; X)$ is the Banach space of all $\mu$-equivalent classes of strongly measurable functions defined on $S$ with values in $X$. In particular, for $X = \mathbb{R}$, $L_p(S; \mathbb{R})$ is denoted by $L_p(S)$. Therefore, the real Banach space $L_p(S)$ is a special case of Bochner spaces.

Let $C$ be a closed and convex subset of this uniformly convex and uniformly smooth Bochner space $L_p(S; X)$. In this paper, we study the following optimization problem: for any given $\varphi \in L_q(S; X^*)$, to find $g \in C$ such that

$$\int_S \langle \varphi(s), g(s) \rangle d\mu(s) = \sup_{f \in C} \int_S \langle \varphi(s), f(s) \rangle d\mu(s). \tag{1.1}$$

Let $J_p \colon L_p(S; X) \to L_q(S; X^*)$ be the normalized duality mapping on $L_p(S; X)$, which will be defined in section 2. Analogously, we consider the following optimization problem: for any $u \in L_p(S; X)$, to find $h \in C$ such that

$$\int_S \langle (J_p u)(s), h(s) \rangle d\mu(s) = \sup_{f \in C} \int_S \langle (J_p u)(s), f(s) \rangle d\mu(s). \tag{1.2}$$

In section 3, we study the optimization problems (1.1) and (1.2) for the general cases that $C$ is a closed and convex subset of the Bochner space $L_p(S; X)$. In section 4, we solve the optimization problems (1.1) and (1.2) for $C$ to be closed and convex cones or closed subspaces in $L_p(S; X)$. In section 5, we consider the optimization problems (1.1) and (1.2) for $C$ to be closed balls in $L_p(S; X)$. For each case, we also study the properties of the inverse images of solution sets.

By definition, every point in the Bochner space $L_p(S; X)$ is the limit of $\mu$-simple functions in $L_p(S; X)$, with respect to the norm $\|\cdot\|_{L_p(S;X)}$. Through sections 3, 4 and 5, we will see that $\mu$-simple functions in $L_p(S; X)$ are the "industry standard" powerful tools for studying optimization problems in Bochner spaces. It is because the importance, in section 2, we recall some useful properties of $\mu$-simple functions in Bochner spaces.

It is well-known that optimization in function spaces is a fundamental object in the theories of optimization and variational analysis. Historically, $L_p(S)$ spaces play the essential role as the main characters in traditional optimization in function spaces (see [9, 11, 15, 21, 26]). Since the standard real Banach space $L_p(S)$ is considered as a special case of Bochner spaces, the traditional optimization in $L_p(S)$ spaces is naturally extended to optimization in Bochner spaces $L_p(S; X)$. This is the motivation of this paper.

The subjects of stochastic optimization, stochastic variational inequality (random variational inequality), stochastic fixed point have been developing rapidly in popular over the last two or three decades (see [3, 8, 19, 22]). In the studying of these topics, the generalized expected vales are multiple applied in the stochastic analysis. Since the generalized expected vales used in these subjects are Bochner integrals, one sees that optimization in Bochner spaces must plays crucial role in stochastic optimization and stochastic variational inequality. Hence, it is indeed very important to study optimization problems in Bochner spaces.

In the techniques of solving optimization problems, variational inequalities, approximation problems and fixed point problems on a closed and convex subset $C$ in a Banach space, the three projections $P_C$ (the standard metric projection), $\pi_C$ (the generalized projection) and $\Pi_C$ (the generalized metric projection) and the normalized duality mapping $J$ play very important roles. In section 2, we recall their definitions and some properties. In particular, we will recall the basic variational principles of these three projections, which are considered as the fundamental theorems in projection theory in Banach spaces.

Since this paper is about optimization in uniformly convex and uniformly smooth Bochner spaces and the closed balls are important underlying subsets, in section 6, we will provide the analytic representations for these projections on closed ball in Bochner spaces.

## 2. Preliminaries

### 2.1. Bochner spaces $(L_p(S; X), \|\cdot\|_{L_p(S;X)})$

In this subsection, we recall the basic definitions and properties of Bochner spaces (see [4−7, 10, 13, 16−17, 20, 23−24]). Let $(S, \mathcal{A}, \mu)$ be a measure space. Throughout this paper, without loss of generality, the measure space $(S, \mathcal{A}, \mu)$ is assumed to be positive and complete. Let $(X, \|\cdot\|_X)$ be a real uniformly convex and uniformly smooth Banach space with topological dual space $(X^*, \|\cdot\|_{X^*})$. Let $\langle \cdot, \cdot \rangle$ denote the real canonical evaluation pairing between the dual space and a given Banach space. For any $A \in \mathcal{A}$ and $x \in X$, let $1_A \otimes x$ denote the $X$-valued simple function on $S$ defined by

$$(1_A \otimes x)(s) = 1_A(s) \otimes x = \begin{cases} x, & \text{if } s \in A, \\ 0, & \text{if } s \notin A, \end{cases} \qquad \text{for any } s \in S,$$

where $1_A$ denotes the characteristic function of $A$ on the space $S$.

For an arbitrary given positive integer $n$, let $\{A_1, A_2, \ldots A_n\}$ be a finite collection of mutually disjoint subsets in $\mathcal{A}$ with $0 < \mu(A_i) < \infty$, for all $i = 1, 2, \ldots n$. Let $\{x_1, x_2, \ldots x_n\} \subseteq X$ and let $\{a_1, a_2, \ldots a_n\}$ be a set of real numbers. Then, $\sum_{i=1}^n a_i(1_{A_i} \otimes x_i)$ is called a $\mu$-simple function from $S$ to $X$ (See Definition 1.1.13 in [10]).

**Remarks 2.1**. 1) In general, a $\mu$-simple function can be written as form $\sum_{i=1}^n a_i(1_{A_i} \otimes x_i)$ with real coefficients $\{a_1, a_2, \ldots a_n\}$.

2) Since the real coefficients $\{a_1, a_2, \ldots a_n\}$ can be included in the points $\{x_1, x_2, \ldots x_n\}$, it follows that a $\mu$-simple function can have the form $\sum_{i=1}^n (1_{A_i} \otimes x_i)$.

Throughout this paper, take positive numbers $p$ and $q$, with $1 < p, q < \infty$ satisfying $\frac{1}{p} + \frac{1}{q} = 1$. Let $(L_p(S; X), \|\cdot\|_{L_p(S;X)})$ be the Lebesgue-Bochner function space called the Bochner space. $L_p(S; X)$ is a real uniformly convex and uniformly smooth Banach space, which is the Banach space of $\mu$-equivalent classes of strongly measurable functions $f: S \to X$ with norm:

$$\|f\|_{L_p(S;X)} = \left( \int_S \|f(s)\|_X^p d\mu(s) \right)^{\frac{1}{p}} < \infty, \text{ for } 1 \leq p < \infty.$$

For $f \in L_p(S; X)$, $f$ takes values in Banach space $X$ as the limit of integrals of simple functions, with respect to the norm $\|\cdot\|_{L_p(S;X)}$. In particular, for $X = \mathbb{R}$, $L_p(S; \mathbb{R})$ is denoted by $L_p(S)$. For easy referee, we list some properties of Bochner integrals and Bochner spaces below. See [6, 10, 17, 24] for more details.

(B₁). $\left\| \int_S f d\mu \right\|_X \leq \int_S \|f\|_X d\mu$, for every $f \in L_p(S, X)$;

(B$_2$). The dual space of $L_p(S; X)$ is $(L_p(S; X))^* = L_q(S; X^*)$ such that, any $\varphi \in L_q(S, X^*)$ defines a bounded linear functional $\varphi \in (L_p(S; X))^*$ by the formula

$$\langle \varphi, f \rangle = \int_S \langle \varphi(\omega), f(\omega) \rangle d\mu(\omega), \text{ for every } f \in L_p(S, X). \tag{2.2}$$

It satisfies

$$\|\varphi\|_{(L_p(S; X))^*} = \|\varphi\|_{L_q(S; X^*)}, \text{ for any } \varphi \in L_q(S, X^*).$$

The geometric properties of Bochner spaces have been studied by many authors, which include the convexity and the smoothness. In particular, the uniformly convexity and uniformly smoothness will be used in this paper.

**Theorem 2.2** (see [17]). Let $(S, \mathcal{A}, \mu)$ be a measure space and let $X$ be a Banach space. For any $p$ with $1 < p < \infty$, we have

(a) $L_P(S; X)$ is uniformly convex, uniformly smooth, if and only if, $X$ is uniformly convex and uniformly smooth;

(b) $L_2(S; X)$ is a Hilbert space $\iff$ $X$ is a Hilbert space.

The uniformly convex and uniformly smooth Bochner space $(L_p(S; X), \|\cdot\|_{L_p(S;X)})$ with dual space $L_q(S; X^*)$ has the following properties and an embedding mapping of $X$ into $L_p(S; X)$, which are studied in [17]. These properties will be multiple used in this paper.

**Proposition 2.3** (see [17]) *For an arbitrary given $A \in \mathcal{A}$ with $0 < \mu(A) < \infty$, for any $x, y \in X$ and for any $1 < p < \infty$, we have*

(i) $\quad \dfrac{1}{\mu(A)^{\frac{1}{p}}}(1_A \otimes x) \in L_p(S; X);$

(ii) $\quad \left\| \dfrac{1}{\mu(A)^{\frac{1}{p}}}(1_A \otimes x) \right\|_{L_p(S;X)} = \|x\|_X;$

(iii) $\quad \left\| \dfrac{1}{\mu(A)^{\frac{1}{p}}}(1_A \otimes x) \pm \dfrac{1}{\mu(A)^{\frac{1}{p}}}(1_A \otimes y) \right\|_{L_p(S;X)} = \|x \pm y\|_X;$

(iv) $\quad$ The mapping $x \to \dfrac{1}{\mu(A)^{\frac{1}{p}}}(1_A \otimes x)$ (isometric) embeds $X$ into $L_p(S; X)$.

## 2.2. The normalized duality mapping in Bochner spaces

Let $X$ be a uniformly convex and uniformly smooth Banach space with dual space $X^*$. The normalized duality mapping $J_X \colon X \to X^*$ is a single-valued mapping satisfying

$$\langle J_X x, x \rangle = \|J_X x\|_{X^*} \|x\|_X = \|x\|_X^2 = \|J_X x\|_{X^*}^2, \text{ for any } x \in X.$$

The normalized duality mapping $J_X$ plays an important role in nonlinear analysis. We list some properties of $J_X$ in uniformly convex and uniformly smooth Banach space for easy reference (see [25] for more details).

(i) $J_X: X \to X^*$ is one to one and onto;
(ii) $J_X$ is a continuous and homogeneous operator and is uniformly continuous on each bounded subset of $X$.

Recall that, in particular, if $X$ is the real Banach space $l_p$ or $L_p(S)$ with $1 < p < \infty$, then the normalized duality mapping $J_X$ in $X$ holds the following analytic representations.

(a) Suppose that $X = l_p$ with $1 < p < \infty$. For any $x = (t_1, t_1, \ldots) \in l_p$ with $x \neq \theta$, we have

$$(J_X x)_n = \frac{|x_n|^{p-1}\text{sign}(x_n)}{\|x\|_{l_p}^{p-2}}, \text{ for } n = 1, 2, \ldots.$$

(b) Suppose that $X = L_p(S)$ with $1 < p < \infty$. For any $h \in L_p(S)$ with $h \neq \theta$, we have

$$(J_X h)(s) = \frac{\|h(s)\|_X^{p-1}\text{sign}(f(s))}{\|f\|_{L_p(S)}^{p-2}}, \text{ for all } s \in S.$$

For the uniformly convex and uniformly smooth Bochner space $L_p(S; X)$ with dual space $L_q(S; X^*)$, throughout this paper, we use the following notations for normalized duality mappings. The normalized duality mappings in $X$ and $L_p(S; X)$ are respectively denoted by $J_X$ and $J_{L_p(S,X)}$. $J_{L_p(S,X)}$ is abbreviated as $J_p$ if there is no confusion caused. Then, by the properties of normalized duality mappings, both $J_X$ and $J_p$ are single valued, one to one and onto continuous mappings. The normalized duality mappings in the dual spaces $X^*$ and $L_q(S; X^*)$ are respectively denoted by $J_{X^*}$ and $J_q^*$, if there is no confusion caused. The following proposition provides the analytic representations and the connections between $J_X$ and $J_p$, which also give an analytic representation of $J_p$.

The normalized duality mappings in the uniformly convex and uniformly smooth Bochner space $L_p(S; X)$ have the following analytic representations and properties, which are proved in [13].

**Proposition 3.1 in [13]**. *For any $f \in L_p(S; X)$ with $f \neq \theta$, we have*

(a) $J_p f \in L_q(S; X^*)$;
(b) $(J_p f)(s) = \dfrac{\|f(s)\|_X^{p-2} J_X(f(s))}{\|f\|_{L_p(S;X)}^{p-2}}$, *for all $s \in S$.*       (2.4)

**Corollary 3.2 in [13]**. *Let $A \in \mathcal{A}$ with $0 < \mu(A) < \infty$. Then, for any $x \in X$ with $x \neq 0$, we have*

$$J_p(1_A \otimes x)(s) = \frac{\mu(A)^{\frac{1}{p}}}{\mu(A)^{\frac{1}{q}}}(1_A \otimes J_X x)(s), \text{ for all } s \in S. \tag{2.5}$$

*It is equivalent to*

$$J_p(\frac{1}{\mu(A)^{\frac{1}{p}}}(1_A \otimes x))(s) = \frac{1}{\mu(A)^{\frac{1}{q}}}(1_A \otimes J_X x)(s), \text{ for all } s \in S. \tag{2.6}$$

**Corollary 3.3 in [13]**. *$J_p$ maps every $\mu$-simple function in $L_p(S; X)$ to $\mu$-simple function in $L_q(S; X^*)$ with respect to the same partition in $S$. Moreover, for an arbitrary given $\mu$-simple function $\sum_{i=1}^n (1_{A_i} \otimes x_i)$ in $L_p(S; X)$, we have*

$$J_p\big(\sum_{i=1}^n (1_{A_i} \otimes x_i)\big)(s) = \frac{1}{\big(\sum_{j=1}^n \|x_j\|_X^p \mu(A_j)\big)^{\frac{1}{q} - \frac{1}{p}}} \sum_{i=1}^n \|x_i\|_X^{p-2} (1_{A_i} \otimes J_X x_i)(s), \text{ for all } s \in S. \quad (2.7)$$

**Corollary 3.4 in [13]**. *For any $f \in L_p(S; X)$, let $\{f_n\}$ be a sequence of $\mu$-simple functions in $L_p(S; X)$ satisfying*

$$f_n \to f, \text{ in } L_p(S; X), \text{ as } n \to \infty.$$

*Then $\{J_p f_n\}$ is a sequence of $\mu$-simple functions in $L_q(S; X^*)$ such that*

$$J_p f_n \to J_p f, \text{ in } L_q(S; X^*), \text{ as } n \to \infty.$$

### 2.3. Three projections in uniformly convex and uniformly smooth Banach Spaces

Let $X$ be a uniformly convex and uniformly smooth Banach space and $C$ a nonempty closed and convex subset of $X$. Let $P_C: X \to C$ denote the (standard) metric projection operator. For $x \in X$, it is defined by

$$\|x - P_C x\|_X \leq \|x - z\|_X, \text{ for all } z \in C.$$

The concepts of generalized projection and generalized metric projection were introduced by Alber in [1] on uniformly convex and uniformly smooth Banach Spaces, which have been extended to general Banach spaces (see [12, 18]).

Let $X$ be a uniformly convex and uniformly smooth Banach space with dual space $X^*$. A Lyapunov functional $V: X^* \times X \to \mathbb{R}_+$ is defined by the following formula:

$$V(\varphi, x) = \|\varphi\|_{X^*}^2 - 2\langle \varphi, x \rangle + \|x\|_X^2, \text{ for any } \varphi \in X^* \text{ and } x \in X.$$

The generalized projection operator $\pi_C: X^* \to C$ is a single-valued mapping defined by

$$V(\varphi, \pi_C \varphi) = \inf_{y \in C} V(\varphi, y), \text{ for } \varphi \in X^*.$$

Based on $\pi_C: X^* \to C$, the generalized metric projection $\Pi_C: B \to C$ is defined by

$$\Pi_C x = \pi_C(J_X x), \text{ for any } x \in X.$$

More prespecify, $\Pi_C: B \to C$ is defined by

$$V(J_X x, \Pi_C x) = \inf_{y \in C} V(J_X x, y), \text{ for any } x \in X.$$

The metric and generalized metric projections have the following connections.

(a) In general, $\Pi_C \neq P_C$;

(b) When $X$ is a Hilbert space, then $\pi_C = \Pi_C = P_C$.

The metric projection $P_C$, the generalized projection $\pi_C$ and the generalized metric projection $\Pi_C$ have many properties. Especially, in uniformly convex and uniformly smooth Banach spaces, they respectively have the following basic variational principles. For $x \in X$, $\varphi \in X^*$ and $y \in C$, we have

$$y = P_C(x) \quad \Longleftrightarrow \quad \langle J_X(x - y), y - z \rangle \geq 0, \text{ for all } z \in C. \tag{2.8}$$

$$y = \pi_C(\varphi) \quad \Longleftrightarrow \quad \langle \varphi - J_X y, y - z \rangle \geq 0, \quad \text{for all } z \in C. \tag{2.9}$$

$$y = \Pi_C(x) \quad \Longleftrightarrow \quad \langle J_X x - J_X y, \ y - z \rangle \geq 0, \text{ for all } z \in C. \tag{2.10}$$

## 3. Optimization problems on closed and convex subsets in Bochner spaces

### 3.1. Optimization problems with respect to dual space

As what is mentioned in the previous section, $L_p(S; X)$ is a uniformly convex and uniformly smooth Bochner space with dual space $L_q(S; X^*)$ with origins $\theta$ and $\theta^*$, respectively, in which $1 < p, q < \infty$ satisfying $\frac{1}{p} + \frac{1}{q} = 1$.

**Definition 3.1**. Let $C$ be a nonempty closed and convex subset in $L_p(S; X)$. For $\varphi \in L_q(S; X^*)$, we consider the following optimization problem with respect to $L_q(S; X^*)$, denoted by OP($C$, $\varphi$): find $g \in C$ such that

$$\int_S \langle \varphi(s), g(s) \rangle d\mu(s) = \sup_{f \in C} \int_S \langle \varphi(s), f(s) \rangle d\mu(s). \tag{3.1}$$

The solution set of OP($C$, $\varphi$) is denoted by $\mathcal{S}_C(\varphi)$. We provide some examples below to demonstrate the solution sets of some optimization problems.

**Example 3.2**. Let $L_p(S; X)$ be a uniformly convex and uniformly smooth Bochner space, in which $(S, \mathcal{A}, \mu)$ is a measure space with $\mu(S) \geq 1$. Take an arbitrary $A \in \mathcal{A}$ with $\mu(A) = 1$. For any given $M > 0$, let

$$C = \left\{ f \in L_p(S; X) : \left( \int_A \|f(s)\|_X^p d\mu(s) \right)^{\frac{1}{p}} \leq M \right\}.$$

Then $C$ is a nonempty closed and convex subset in $L_p(S; X)$. Take an arbitrary $x^* \in X^*$ with $\|x^*\|_{X^*} = 1$. It follows that $1_A \otimes x^* \in L_q(S; X^*)$. Then the solution set $\mathcal{S}_C(1_A \otimes x^*)$ of (3.1) is nonempty and it is

$$\mathcal{S}_C(1_A \otimes x^*) = \left\{ g \in C : \int_A \langle x^*, g(s) \rangle d\mu(s) = M \right\}. \tag{3.2}$$

*Proof*. For any $g \in C$, (3.1) could be rewritten as,

$$\langle \varphi, g \rangle = \sup_{f \in C} \langle \varphi, f \rangle.$$

If $\int_A \langle x^*, g(s) \rangle d\mu(s) = M$, then

$$\begin{aligned}
\langle 1_A \otimes x^*, g \rangle &= \int_S \langle (1_A \otimes x^*)(s), g(s) \rangle d\mu(s) \\
&= \int_A \langle x^*, g(s) \rangle d\mu(s) \\
&= M.
\end{aligned} \tag{3.3}$$

For any $f \in C$, by $\mu(A) = 1$ and $\|x^*\|_{X^*} = 1$, we have

$$\begin{aligned}
\langle 1_A \otimes x^*, f \rangle &= \int_S \langle (1_A \otimes x^*)(s), f(s) \rangle d\mu(s) \\
&= \int_A \langle x^*, f(s) \rangle d\mu(s) \\
&\leq \int_A \|x^*\|_{X^*} \|f(s)\|_X d\mu(s) \\
&\leq \left( \int_A \|f(s)\|_X^p d\mu(s) \right)^{\frac{1}{p}} \\
&\leq M.
\end{aligned} \tag{3.4}$$

Then, For any $g \in C$, (3.3) and (3.4) imply that

$$\langle 1_A \otimes x^*, g \rangle = \int_A \langle x^*, g(s) \rangle d\mu(s) = M \implies g \in \mathcal{S}_C(1_A \otimes x^*).$$

This implies

$$\left\{ g \in C : \int_A \langle x^*, g(s) \rangle d\mu(s) = M \right\} \subseteq \mathcal{S}_C(1_A \otimes x^*). \tag{3.5}$$

On the other hand, let $x = J^* x^*$, in which $x^* \in X^*$ is given in this example with $\|x^*\|_{X^*} = 1$. Then $\|x\|_X = \|x^*\|_{X^*} = 1$. Define a function $h \in L_p(S; X)$ by

$$h(s) = M(1_A \otimes x)(s), \text{ for all } s \in S.$$

By $\mu(A) = 1$ and $\|x\|_X = 1$, we calculate $\|h\|_{L_p(S;X)} = M$. This implies $h \in C$. Then, we have

$$\begin{aligned}
\langle 1_A \otimes x^*, h \rangle &= \int_S \langle (1_A \otimes x^*)(s), M(1_A \otimes x)(s) \rangle d\mu(s) \\
&= \int_A M \langle x^*, x \rangle d\mu(s) \\
&= M.
\end{aligned} \tag{3.6}$$

By (3.5), this implies that $h \in \mathcal{S}_C(1_A \otimes x^*)$ and the optimization (3.1) attains the maximum. Then, for any $u \in \mathcal{S}_C(1_A \otimes x^*) \subseteq C$, by (3.6), we must have

$$M \geq \|u\|_{L_p(S;X)} \geq \langle 1_A \otimes x^*, u \rangle \geq \langle 1_A \otimes x^*, h \rangle = M.$$

It follows that $\langle 1_A \otimes x^*, u \rangle = M$. This implies

$$\mathcal{S}_C(1_A \otimes x^*) \subseteq \Big\{ g \in C \colon \int_A \langle x^*, g(s) \rangle d\mu(s) = M \Big\}. \tag{3.7}$$

By (3.5) and (3.7), it proves (3.2). Since $h \in \mathcal{S}_C(1_A \otimes x^*)$, which shows that

$$\mathcal{S}_C(1_A \otimes x^*) \neq \emptyset. \qquad \square$$

**Example 3.3**. Let $L_p(S; X)$ be a uniformly convex and uniformly smooth Bochner space, in which $(S, \mathcal{A}, \mu)$ is a measure space with $\mu(S) \geq 3$ (this is just for the simpler calculation) and $(X, \|\cdot\|_X)$ is a uniformly convex and uniformly smooth Banach space with dimension greater than or equal to 3. Take arbitrary mutually disjoint elements $A_1, A_2, A_3$ in $\mathcal{A}$ with $\mu(A_i) = 1$, and take arbitrary three linearly independent points $x_1, x_2, x_3 \in X$ with $\|x_i\|_X = 1$ and $x_i^* = J_X(x_i)$, which satisfies $\|x_i^*\|_{X^*} = \|x_i\|_X = 1$, for $i = 1, 2, 3$.

Let $C$ be the (closed) convex hell of $\{(1_{A_1} \otimes x_1), (1_{A_2} \otimes x_2), (1_{A_3} \otimes x_3)\}$

$$C = \{ t_1(1_{A_1} \otimes x_1) + t_2(1_{A_2} \otimes x_2) + t_3(1_{A_3} \otimes x_3) \in X \colon t_i \geq 0, \text{ for } i = 1, 2, 3 \text{ with } t_1 + t_2 + t_3 = 1 \}.$$

Then $C$ is a nonempty closed and convex subset of $L_p(S; X)$. Take four points in $L_q(S; X^*)$ as

$$\varphi = (1_{A_1} \otimes x_1^*) + (1_{A_2} \otimes x_2^*) + (1_{A_3} \otimes x_3^*),$$
$$\psi = (1_{A_1} \otimes x_1^*) + (1_{A_2} \otimes x_2^*),$$
$$\gamma = (1_{A_1} \otimes x_1^*),$$
$$\lambda = (1_{A_1} \otimes x_1^*) - (1_{A_2} \otimes x_2^*) + (1_{A_3} \otimes x_3^*).$$

Then, we have

(i)  $\mathcal{S}_C(\varphi) = C$ such that
$$\langle \varphi, f \rangle = 1, \text{ for every } f \in C;$$

(ii)  $\mathcal{S}_C(\psi) = \{ t_1(1_{A_1} \otimes x_1) + t_2(1_{A_2} \otimes x_2) \in C \colon t_i \geq 0, \text{ for } i = 1, 2 \text{ with } t_1 + t_2 = 1 \}$. Then
$$\langle \psi, t_1(1_{A_1} \otimes x_1) + t_2(1_{A_2} \otimes x_2) \rangle = 1, \text{ for any } t_i \geq 0, \text{ for } i = 1, 2 \text{ with } t_1 + t_2 = 1;$$

(iii)  $\mathcal{S}_C(\gamma) = \{ 1_{A_1} \otimes x_1 \};$

(iv)  $\mathcal{S}_C(\lambda) = \{ t_1(1_{A_1} \otimes x_1) + t_3(1_{A_3} \otimes x_3) \in C \colon t_i \geq 0, \text{ for } i = 1, 3 \text{ with } t_1 + t_3 = 1 \}.$

*Proof*. The proof of this example is similar to the proof of Example 3.2 and omitted here. $\square$

**Example 3.4**. Let $L_p(S; X)$ be a uniformly convex and uniformly smooth Bochner space, in which $(X, \|\cdot\|_X)$ is a uniformly convex and uniformly smooth Banach space with dimension greater than or equal to 3. Take mutually disjoint elements $A_1, A_2, A_3$ in $\mathcal{A}$ with $\mu(A_i) > 0$, and take three linearly independent points $x_1, x_2, x_3 \in X$ with $\|x_i\|_X = 1$, for $i = 1, 2, 3$. Define $u \in L_p(S; X)$ by

$$u = \frac{25}{\mu(A_1)^{\frac{1}{p}}} \left(1_{A_1} \otimes x_1\right) + \frac{37}{\mu(A_2)^{\frac{1}{p}}} \left(1_{A_2} \otimes x_2\right) + \frac{77}{\mu(A_3)^{\frac{1}{p}}} \left(1_{A_3} \otimes x_3\right)$$

Let

$$K = \{tu \in X : t \geq 0\}.$$

Then $K$ is a ray in $L_p(S; X)$ with ending point at $\theta$ and toward with direction $u$, which is a closed, convex and pointed cone in $L_p(S; X)$. Let $x_i^* = J_X(x_i)$ with $\|x_i^*\|_{X^*} = \|x_i\|_X = 1$, for $i = 1, 2, 3$. Then, define $\varphi$, $\psi$, $\gamma \in L_q(S; X^*)$ by

$$\Phi = \frac{-9}{\mu(A_1)} \left(1_{A_1} \otimes x_1^*\right) + \frac{4}{\mu(A_2)} \left(1_{A_2} \otimes x_2^*\right) + \frac{1}{\mu(A_3)} \left(1_{A_3} \otimes x_3^*\right),$$

$$\psi = \frac{-9}{\mu(A_1)} \left(1_{A_1} \otimes x_1^*\right) + \frac{-1}{\mu(A_3)} \left(1_{A_3} \otimes x_3^*\right),$$

$$\gamma = \frac{9}{\mu(A_1)} \left(1_{A_1} \otimes x_1^*\right) + \frac{4}{\mu(A_2)} \left(1_{A_2} \otimes x_2^*\right).$$

Then, we have

(i)     $\mathcal{S}_K(\Phi) = K$;
(ii)    $\mathcal{S}_K(\psi) = \{\theta\}$;
(iii)   $\mathcal{S}_K(\gamma) = \emptyset$.

*Proof.* Proof of (i). For any $tu \in K$ with $t \geq 0$, we calculate

$$\langle \Phi, tu \rangle$$

$$= t\left(\int_{A_1} \frac{(-9)25}{\mu(A_1)} \langle x_1^*, x_1 \rangle d\mu(s) + \int_{A_2} \frac{(4)37}{\mu(A_2)} \langle x_2^*, x_2 \rangle d\mu(s) + \int_{A_3} \frac{77}{\mu(A_3)} \langle x_3^*, x_3 \rangle d\mu(s)\right)$$

$$= t\left((-9)25 + (4)37 + 77\right)$$
$$= 0.$$

This reduces that $tu \in \mathcal{S}_K(\Phi)$, for any $tu \in K$ with $t \geq 0$.

Proofs of (ii) and (iii). Similar to the proof of (i), for any $tu \in K$, we have

$$\langle \psi, tu \rangle = t\left((-9)25 - 77\right) < 0, \text{ for any } t > 0,$$

and

$$\langle \gamma, tu \rangle = t\left((9)25 + (4)37\right) \to \infty, \text{ as } t \to \infty.$$

(ii) and (iii) follow immediately.                                   □

From the above three examples, $OP(C, \varphi)$ has the following properties.

**Proposition 3.5.** *Let $C$ be a nonempty closed and convex subset of $L_p(S; X)$. Then, we have*

(a) $\mathcal{S}_C(\theta^*) = C$;
(b) *For* $\varphi \in L_q(S; X^*)$, *if* $\mathcal{S}_C(\varphi) \neq \emptyset$, *then* $\mathcal{S}_C(\varphi)$ *is a closed and convex subset of C.*

*Proof.* The proof of this proposition is straight forward and omitted here. $\square$

**Proposition 3.6.** *Let C be a nonempty closed, convex and bounded subset of* $L_p(S; X)$. *Then, for any* $\varphi \in (L_p(S; X))^*$, $\mathcal{S}_C(\varphi)$ *is a nonempty closed and convex subset of C.*

*Proof.* Since $C$ is a nonempty closed, convex and bounded subset in the uniformly convex and uniformly smooth Bochner space $L_p(S; X)$, then it is weakly compact. For any $\varphi \in (L_p(S; X))^*$, the functional

$$\langle \varphi, g \rangle = \int_S \langle \varphi(s), g(s) \rangle d\mu(s), \text{ for any } g \in L_p(S; X)$$

is continuous; and therefore, it attains its maximum value on $C$. This implies $\mathcal{S}_C(\varphi) \neq \emptyset$. The proof of the closeness and convexity of $\mathcal{S}_C(\varphi)$ is straight forward and omitted here. $\square$

**Definition 3.7.** Let $C$ be a nonempty closed and convex subset of $L_p(S; X)$. Let

$$\mathcal{S}_C^{-1}(g) = \left\{ \varphi \in L_q(S; X^*) \colon \int_S \langle \varphi(s), g(s) \rangle d\mu(s) = \sup_{f \in C} \int_S \langle \varphi(s), f(s) \rangle d\mu(s) \right\}.$$

Notice that, for $g \in C$, $\mathcal{S}_C^{-1}(g)$ could be rewritten as

$$\mathcal{S}_C^{-1}(g) = \left\{ \varphi \in L_q(S; X^*) \colon \langle \varphi, g \rangle = \sup_{f \in C} \langle \varphi, f \rangle \right\} = \left\{ \varphi \in L_q(S; X^*) \colon g \in \mathcal{S}_C(\varphi) \right\}.$$

**Proposition 3.8.** *Let C be a nonempty closed and convex subset of* $L_p(S; X)$. *Then, we have*

(a) $\theta^* \in \mathcal{S}_C^{-1}(g)$, *for every* $g \in C$.
(b) *For* $g \in C$, *if* $\mathcal{S}_C^{-1}(g) \supsetneq \{\theta^*\}$, *then* $\mathcal{S}_C^{-1}(g)$ *is a closed and convex cone in* $L_q(S; X^*)$ *with vertex at* $\theta^*$.

*Proof.* The proof of this proposition is straight forward and omitted here. $\square$

**Definition 3.9.** For any given $g \in C$, if $\mathcal{S}_C^{-1}(g) \supsetneq \{\theta^*\}$, then $g$ is said to be optimal. Otherwise, $g$ is said to be none-optimal. The set of all optimal points in $C$ is denoted by $O(C)$; and the set of all none-optimal points in $C$ is denoted by $N(C)$.

By Proposition 3.8, we have that $\{O(C), N(C)\}$ forms a partition of $C$, that is,

$$O(C) \cup N(C) = C \quad \text{and} \quad O(C) \cap N(C) = \emptyset.$$

### 3.2. Optimization problems with respect to the underlying space

Let $L_p(S; X)$ be a uniformly convex and uniformly smooth Bochner space with dual space $L_q(S; X^*)$. It is known that the normalized duality mapping $J_p: L_p(S; X) \to L_q(S; X^*)$ is a one to one, onto and continuous mapping. Analogously to Definition 3.1, we consider the following optimization problem with respect to $L_p(S; X)$.

**Definition 3.10**. Let $L_p(S; X)$ be a uniformly convex and uniformly smooth Bochner space with dual space $L_q(S; X^*)$ and let $C$ be a nonempty closed and convex subset of $L_p(S; X)$. We define an optimization problem with respect to $L_p(S; X)$ as follows: for any $u \in L_p(S; X)$, to find $h \in C$ such that

$$\int_S \langle (J_p u)(s), h(s) \rangle d\mu(s) = \sup_{f \in C} \int_S \langle (J_p u), f(s) \rangle d\mu(s). \tag{3.8}$$

For $g \in C$, let

$$\mathcal{S}_C^{-*}(g) = \left\{ h \in L_p(S; X): \int_S \langle (J_p h)(s), g(s) \rangle d\mu(s) = \sup_{f \in C} \int_S \langle (J_p h)(s), f(s) \rangle d\mu(s) \right\}$$
$$= \left\{ h \in L_p(S; X): g \in \mathcal{S}_C(J_p h) \right\}.$$

From Definitions 3.7 and 3.10, $\mathcal{S}_C^{-1}$ is a set-valued mapping from $C$ to $L_q(S; X^*)$ and $\mathcal{S}_C^{-*}$ is a set-valued mapping from $C$ to $L_p(S; X)$, which can be rewritten as

$$\mathcal{S}_C^{-*}(g) = \left\{ h \in L_p(S; X): g \in \mathcal{S}_C(J_p h) \right\}.$$

$\mathcal{S}_C^{-1}$ and $\mathcal{S}_C^{-*}(g)$ are inverse images of the solution sets of the optimization problems (3.1) and (3.8), respectively. They have the following connections.

**Lemma 3.11**. *Let $C$ be a nonempty closed and convex subset of $L_p(S; X)$. Then, for $g \in C$, we have*

$$g \in \mathcal{S}_C^{-*}(g) \iff \|g\|_{L_p(S;X)}^2 = \sup_{f \in C} \int_S \langle (J_p g), f(s) \rangle d\mu(s).$$

*Proof.* By Definition (3.8), we have

$$g \in \mathcal{S}_C^{-*}(g) \iff \int_S \langle (J_p g)(s), g(s) \rangle d\mu(s) = \sup_{f \in C} \int_S \langle (J_p g), f(s) \rangle d\mu(s).$$

Then, this lemma follows from

$$\|g\|_{L_p(S;X)}^2 = \int_S \langle (J_p g)(s), g(s) \rangle d\mu(s). \qquad \square$$

**Lemma 3.12**. *Let $C$ be a nonempty closed and convex subset of $L_p(S; X)$. Then, for any $g \in C$, we have*

$$\mathcal{S}_C^{-*}(g) = J_q^*(\mathcal{S}_C^{-1}(g)) \quad \text{and} \quad \mathcal{S}_C^{-1}(g) = J_p(\mathcal{S}_C^{-*}(g)).$$

*Proof.* Since $L_p(S; X)$ and $L_q(S; X^*)$ both are uniformly convex and uniformly smooth Banach spaces, then $J_p$ and $J_q^*$ both are 1-1 and onto mappings such that

$$J_q^* J_p = I_{L_p(S;X)} \quad \text{and} \quad J_p J_q^* = I_{L_q(S;X^*)}.$$

This lemma follows immediately.                                                    □

**Theorem 3.13**. *Let C be a nonempty closed and convex subset of $L_p(S; X)$. Then, for any $g \in C$, we have*

(a) $\Theta \in \mathcal{S}_C^{-*}(g)$;

(b) *If $\mathcal{S}_C^{-*}(g) \supsetneq \{\Theta\}$, then $\mathcal{S}_C^{-*}(g)$ is a closed cone with vertex at $\Theta$ in $L_p(S; X)$;*

(c) *In general, $\mathcal{S}_C^{-*}(g)$ is not convex.*

*Proof*. Since $L_p(S; X)$ and $L_q(S; X^*)$ both are uniformly convex and uniformly smooth Banach spaces, then $J_p$ and $J_q^*$ both are 1-1, onto and continuous mappings. By Proposition 3.8, the positive homogenous of $J_q^*$ and from Proposition 3.11, $\mathcal{S}_C^{-*}(g) = J_q^*(\mathcal{S}_C^{-1}(g))$, parts (a, b) follow immediately. So, we only prove part (c).

We prove part (c) by construct a counter example, which is similar to Example 3.3 for $p = 3$. Let $L_3(S; X)$ be a uniformly convex and uniformly smooth Bochner space, in which $(S, \mathcal{A}, \mu)$ is a measure space with $\mu(S) \geq 3$ and $(X, \|\cdot\|_X)$ is a uniformly convex and uniformly smooth Banach space with dimension greater than or equal to 3. Take arbitrary mutually disjoint elements $A_1$, $A_2$, $A_3$ in $\mathcal{A}$ with $\mu(A_i) = 1$, and arbitrary three linearly independent points $x_1, x_2, x_3 \in X$ with $\|x_i\|_X = 1$ and $x_i^* = J_X(x_i)$ such that $\|x_i^*\|_{X^*} = \|x_i\|_X = 1$, for $i = 1, 2, 3$. Let

$$g = 25(1_{A_1} \otimes x_1) + 37(1_{A_2} \otimes x_2) + 77(1_{A_3} \otimes x_3) \in L_3(S; X).$$

Then $g$ is a $\mu$-simple function in $L_3(S; X)$. Let $C$ be the closed segment in $L_3(S; X)$ with ending points at $\Theta$ and $g$:

$$C = \{tg \in L_3(S; X) : t \in [0, 1]\}.$$

Then $C$ is a nonempty closed and convex subset of $L_3(S; X)$. Take $u, v \in L_3(S; X)$ as follows

$$u = 3(1_{A_1} \otimes x_1) - 2(1_{A_2} \otimes x_2) - (1_{A_3} \otimes x_3),$$

and

$$v = (1_{A_1} \otimes x_1) - 3(1_{A_2} \otimes x_2) + 2(1_{A_3} \otimes x_3).$$

By $\mu(A_i) = 1$, $x_i^* = J_X(x_i)$ with $\|x_i^*\|_{X^*} = \|x_i\|_X = 1$, for $i = 1, 2, 3$, and by (2.7), we have

$$J_3(u) = \frac{9}{\sqrt[3]{36}}\left(1_{A_1} \otimes x_1^*\right) - \frac{4}{\sqrt[3]{36}}\left(1_{A_2} \otimes x_2^*\right) - \frac{1}{\sqrt[3]{36}}(1_{A_3} \otimes x_3^*) \in L_{\frac{3}{2}}(S; X^*),$$

and

$$J_3(v) = \frac{1}{\sqrt[3]{36}}\left(1_{A_1} \otimes x_1^*\right) - \frac{9}{\sqrt[3]{36}}\left(1_{A_2} \otimes x_2^*\right) + \frac{4}{\sqrt[3]{36}}(1_{A_3} \otimes x_3^*) \in L_{\frac{3}{2}}(S; X^*).$$

Then, we calculate

$$\langle J_3(u), g - tg \rangle$$
$$= \int_S \langle (J_3 u)(s), (g(s) - tg(s)) d\mu(s)$$

$$= (1-t)\int_S \left\langle \left(\tfrac{9}{\sqrt[3]{36}}\left(1_{A_1}\otimes x_1^*\right) - \tfrac{4}{\sqrt[3]{36}}\left(1_{A_2}\otimes x_2^*\right) - \tfrac{1}{\sqrt[3]{36}}\left(1_{A_3}\otimes x_3^*\right)\right)(s),\ \left(25(1_{A_1}\otimes x_1) + 37(1_{A_2}\otimes x_2) + 77(1_{A_3}\otimes x_3)\right)(s)\right\rangle d\mu(s)$$

$$= (1-t)\left\langle \left(\tfrac{9}{\sqrt[3]{36}},\ \tfrac{-4}{\sqrt[3]{36}},\ \tfrac{-1}{\sqrt[3]{36}}\right),\ (25, 37, 77)\right\rangle$$

$$= 0,\ \text{for every } tg \in C \text{ with } t \in [0, 1].$$

By Definition 3.10, this implies $u \in \mathcal{S}_C^{-*}(g)$. We can Similarly to calculate $v \in \mathcal{S}_C^{-*}(g)$. Hence, we obtain that both $u$ and $v$ are in $\mathcal{S}_C^{-*}(g)$. Now, we take a convex combination $h$ of $u$ and $v$ as:

$$h = \tfrac{2}{3}u + \tfrac{1}{3}v = \tfrac{7}{3}(1_{A_1}\otimes x_1) - \tfrac{7}{3}(1_{A_2}\otimes x_2) - 0(1_{A_3}\otimes x_3) \in L_3(S; X).$$

By (2.7), we calculate

$$J_3(h) = \frac{\left|\tfrac{7}{3}\right|^2}{\tfrac{7}{3}\sqrt[3]{2}}\left(1_{A_1}\otimes x_1^*\right) - \frac{\left|\tfrac{7}{3}\right|^2}{\tfrac{7}{3}\sqrt[3]{2}}\left(1_{A_2}\otimes x_2^*\right) - \frac{0}{\tfrac{7}{3}\sqrt[3]{2}}\left(1_{A_3}\otimes x_3^*\right)$$
$$= \frac{7\sqrt[3]{4}}{6}\left[\left(1_{A_1}\otimes x_1^*\right) - \left(1_{A_2}\otimes x_2^*\right)\right] \in L_{\frac{3}{2}}(S; X^*).$$

This implies

$$\langle J_3(h), g - tg\rangle$$

$$= \int_S \langle (J_3 h)(s), (g(s) - tg(s))d\mu(s)$$

$$= (1-t)\int_S \left\langle \tfrac{7\sqrt[3]{4}}{6}\left(\left(1_{A_1}\otimes x_1^*\right) - \left(1_{A_2}\otimes x_2^*\right)\right)(s),\ \left(25(1_{A_1}\otimes x_1) + 37(1_{A_2}\otimes x_2) + 77(1_{A_3}\otimes x_3)\right)(s)\right\rangle d\mu(s)$$

$$= (1-t)\tfrac{7\sqrt[3]{4}}{6}\langle (1, -1, 0), (25, 37, 77)\rangle$$

$$= -14\sqrt[3]{4}\,(1-t)$$

$$< 0,\ \text{for every } tg \in C \text{ with } t \in (0, 1).$$

This shows that $h \notin \mathcal{S}_C^{-*}(g)$. Hence $\mathcal{S}_C^{-*}(g)$ is not convex, which proves this proposition. $\quad\square$

## 4. Optimization problems on cones and subspaces in Bochner spaces

As what is mentioned in the previous sections, throughout this section, $L_p(S; X)$ is a uniformly convex and uniformly smooth Bochner space with dual space $L_q(S; X^*)$ with origins $\theta$ and $\theta^*$, respectively, in which $1 < p, q < \infty$ satisfying $\frac{1}{p} + \frac{1}{q} = 1$. Let $K$ be a nonempty closed and convex cone in $L_p(S; X)$. Then, for any $g \in K$, by Proposition 3.8, either $\mathcal{S}_K^{-1}(g) = \{\theta^*\}$, or, $\mathcal{S}_K^{-1}(g)$ is a closed and convex cone in $L_q(S; X^*)$ with vertex at $\theta^*$. For any $\varphi \in L_q(S; X^*)$ and $f \in L_p(S; X)$, we write

$$\varphi \perp f \quad \text{if and only if} \quad \langle \varphi, f\rangle = 0.$$

Let $A$ be a nonempty subset in $L_p(S; X)$, we write

$$\varphi \perp A \quad \text{if and only if} \quad \varphi \perp f, \text{ for every } f \in A.$$

**Theorem 4.1**. *Let $K$ be a closed and convex cone in $L_p(S; X)$ with vertex at $v \in L_p(S; X)$. For any $\varphi \in L_q(S; X^*)$, we have*

(a) $g \in \mathcal{S}_K(\varphi) \implies \varphi \perp (g - v)$; *that is,*

$$\langle \varphi, g \rangle = \langle \varphi, v \rangle, \text{ for any } g \in \mathcal{S}_K(\varphi);$$

(b) $\mathcal{S}_K(\varphi) \neq \emptyset \implies v \in \mathcal{S}_K(\varphi)$;

(c) *If there is $g \in \mathcal{S}_K(\varphi)$ with $g \neq v$, then $\mathcal{S}_K(\varphi)$ is a closed and convex cone in $L_p(S; X)$ contained in $K$ with the common vertex at $v \in L_p(S; X)$;*

(d) $\mathcal{S}_K(\varphi) \neq \emptyset \implies \varphi \perp (\mathcal{S}_K(\varphi) - v)$; *In particular, if $v = \theta$, then*

$$\varphi \perp \mathcal{S}_K(\varphi), \text{ for any } \varphi \in L_q(S; X^*) \text{ with } \mathcal{S}_K(\varphi) \neq \emptyset.$$

(e) $\cap\{\mathcal{S}_K(\psi): \psi \in L_q(S; X^*) \text{ with } \mathcal{S}_K(\psi) \neq \emptyset\} \neq \emptyset$, *and*

$$v \in \cap\{\mathcal{S}_K(\psi): \psi \in L_q(S; X^*) \text{ with } \mathcal{S}_K(\psi) \neq \emptyset\}.$$

*Proof.* (a) It is trivial for $g = v$. So we assume that $g \neq v$. By $g \in \mathcal{S}_K(\varphi)$ with $g \neq v$, we have

$$\langle \varphi, g - (v + t(g - v)) \rangle \geq 0, \text{ for all } t \geq 0.$$

It follows

$$(1-t)\langle \varphi, g - v \rangle \geq 0, \text{ for all } t \geq 0.$$

Since $g - v \neq \theta$, this implies

$$\langle \varphi, g - v \rangle = 0. \tag{4.1}$$

(b) If $\mathcal{S}_K(\varphi) \neq \emptyset$, we take $g \in \mathcal{S}_K(\varphi)$, by (4.1), for every $f \in K$, we have

$$\begin{aligned}
0 &\leq \langle \varphi, g - f \rangle \\
&= \langle \varphi, g - v + v - f \rangle \\
&= \langle \varphi, v - f \rangle.
\end{aligned}$$

This implies that

$$v \in \mathcal{S}_K(\varphi), \text{ for any } \varphi \in L_q(S; X^*) \text{ with } \mathcal{S}_K(\varphi) \neq \emptyset. \tag{4.2}$$

(c). For any $h \in \mathcal{S}_K(\varphi)$ with $h \neq v$, by (a), $h$ must satisfy (4.1). Then, by (b), $v \in \mathcal{S}_K(\varphi)$, for all $t \geq 0$, we have

$$\begin{aligned}
&\langle \varphi, (v + t(h - v)) - f \rangle \\
&= t\langle \varphi, h - v \rangle + \langle \varphi, v - f \rangle \\
&= \langle \varphi, v - f \rangle \\
&\geq 0, \text{ for all } f \in K.
\end{aligned}$$

This implies

$$v + t(h - v) \in \mathcal{S}_K(\varphi), \text{ for all } t \geq 0.$$

It follows that $\mathcal{S}_K(\varphi)$ is a cone contained in $K$ with the common vertex at $v \in L_p(S; X)$. The closeness and convexity of $\mathcal{S}_K(\varphi)$ follow from part (b) of Proposition 3.5. (d) follows from (a).

(e). By (4.2), this implies that

$$v \in \cap\{\mathcal{S}_K(\psi): \psi \in L_q(S; X^*) \text{ with } \mathcal{S}_K(\psi) \neq \emptyset\}.$$

By (a) of Proposition 3.5, $\mathcal{S}_K(\theta^*) = K$ and by (4.2) again, it implies

$$\cap\{\mathcal{S}_K(\psi): \psi \in L_q(S; X^*) \text{ with } \mathcal{S}_K(\psi) \neq \emptyset\} \neq \emptyset. \qquad \square$$

**Corollary 4.2**. *Let $K$ be a closed and convex cone in $L_p(S; X)$ with vertex at $v \in L_p(S; X)$. Then, for any $\varphi \in L_q(S; X^*)$,*

$$\mathcal{S}_K(\varphi) \neq \emptyset \qquad \Longleftrightarrow \qquad v \in \mathcal{S}_K(\varphi).$$

*Proof*. From part (b) of Theorem 4.1, $\mathcal{S}_K(\varphi) \neq \emptyset \Longrightarrow v \in \mathcal{S}_K(\varphi)$. The other direction is evident.$\square$

**Corollary 4.3**. *Let $K$ be a closed and convex cone in $L_p(S; X)$ with vertex at $\theta$. Then, for any $\varphi \in L_q(S; X^*)$,*

$$\mathcal{S}_K(\varphi) \neq \emptyset \qquad \Longrightarrow \qquad \varphi \perp \mathcal{S}_K(\varphi).$$

*Proof*. This is an immediate consequence of part (d) of Theorem 4.1. $\qquad \square$

As a corollary of Corollary 4.3, we have the following results immediately.

**Corollary 4.4**. *Let $K$ be a closed and convex cone in $L_p(S; X)$ with vertex at $\theta$. Then,*

$$J_q^* \varphi \notin \mathcal{S}_K(\varphi), \text{ for any } \varphi \in L_q(S; X^*) \text{ with } \varphi \neq \theta.$$

*Proof*. If $J_q^* \varphi \notin K$, it is clear to have $J_q^* \varphi \notin \mathcal{S}_K(\varphi)$. If $J_q^* \varphi \in K$, and $\mathcal{S}_K(\varphi) = \emptyset$, this corollary is evident. If $J_q^* \varphi \in K$, and $\mathcal{S}_K(\varphi) \neq \emptyset$, by Corollary 4.3, $\varphi \perp \mathcal{S}_K(\varphi)$. Since $\langle \varphi, J_q^* \varphi \rangle = \|\varphi\|_{X^*}^2 > 0$, this implies $J_q^* \varphi \notin \mathcal{S}_K(\varphi)$. $\qquad \square$

**Corollary 4.5**. *Let $K$ be a closed and convex cone in $L_p(S; X)$ with vertex at $\theta \in L_p(S; X)$. Let $\varphi = \sum_{i=1}^n (1_{A_i} \otimes x_i^*)$ be an $\mu$-simple functional in $L_q(S; X^*)$, in which $\{A_1, A_2, \dots A_n\}$ is a finite collection of mutually disjoint subsets in $\mathcal{A}$ with $0 < \mu(A_i) < \infty$ and $x_i^* \in X^*$, for $i = 1, \dots n$.*

(i)      *If $v \neq \theta$, then,*

$$\mathcal{S}_K(\varphi) \neq \emptyset \Longleftrightarrow \sum_{i=1}^n \int_{A_i} \langle x_i^*, v(s) \rangle d\mu(s) \geq \sum_{i=1}^n \int_{A_i} \langle x_i^*, f(s) \rangle d\mu(s), \text{ for all } f \in K. \quad (4.3)$$

(ii)      *If $v = \theta$, then, $\mathcal{S}_K(\varphi) \neq \emptyset$, if and only if*

$$\mathcal{S}_K(\varphi) \neq \emptyset \qquad \Longleftrightarrow \qquad \sum_{i=1}^n \int_{A_i} \langle x_i^*, f(s) \rangle d\mu(s) \leq 0, \text{ for all } f \in K.$$

*Proof.* Proof of (i). For the given $\varphi = \sum_{i=1}^{n} (1_{A_i} \otimes x_i^*) \in L_q(S; X^*)$, condition (4.3) is equivalent to $v \in \mathcal{S}_K(\varphi)$. Then (i) of this corollary follows from Corollary 4.2 immediately. Part (ii) is a special case of (i) for $v = \theta$. □

**Theorem 4.6.** *Let $D$ be a proper closed subspace of $L_p(S; X)$. Then, for any $\varphi \in L_q(S; X^*)$,*

$$\mathcal{S}_D(\varphi) \neq \emptyset \quad \Longleftrightarrow \quad \mathcal{S}_D(\varphi) = D \quad \Longleftrightarrow \quad \varphi \perp D.$$

*Proof.* Since $D$ is a proper closed subspace of $L_p(S; X)$, for any point $v \in D$, $D$ is a closed and convex cone with vertex at $v$. By part (b) in Theorem 4.1,

$$\mathcal{S}_D(\varphi) \neq \emptyset \quad \Longrightarrow \quad v \in \mathcal{S}_D(\varphi), \text{ for any point } v \in D.$$

By Corollary 4.3, this implies

$$\mathcal{S}_D(\varphi) \neq \emptyset \quad \Longrightarrow \quad \mathcal{S}_D(\varphi) = D \Longrightarrow \varphi \perp D.$$

The direction "⟸" is clear. □

For any $\varphi \in L_q(S; X^*)$, we write

$$\mathrm{SP}(\varphi) = \{s \in S \colon \varphi(s) \neq \text{the origin in } X^*\}.$$

$\mathrm{SP}(\varphi)$ is called the support of $\varphi$.

**Example 4.7.** Let $\{A_1, A_2, \dots A_n\}$ be a finite collection of mutually disjoint subsets in $\mathcal{A}$ with $0 < \mu(A_i) < \infty$, for every $i = 1, 2, \dots n$. Let $\{x_1, x_2, \dots x_n\}$ be a finite subset in $X$. Let $D$ be the closed subspace of $L_p(S; X)$ spanned by $\{1_{A_1} \otimes x_1, 1_{A_2} \otimes x_2, ..., 1_{A_n} \otimes x_n\}$. Then, for any $\varphi \in L_q(S; X^*)$,

$$\mathrm{SP}(\varphi) \cap (\cup_{i=1}^n A_i) = \emptyset \quad \Longrightarrow \quad \mathcal{S}_D(\varphi) = D.$$

*Proof.* For any $f \in D$, $f = \sum_{i=1}^{n} t_i (1_{A_n} \otimes x_n)$, for real numbers $t_i$, for $i = 1, 2, \dots n$. It follows

$$\langle \varphi, f \rangle = \int_S \langle \varphi(s), f(s) \rangle d\mu(s) = \sum_{i=1}^{n} \int_{A_i} \langle \varphi(s), x_i \rangle d\mu(s) = 0, \text{ for any } f \in D.$$

This implies $\varphi \perp D$. By Theorem 4.6, we obtain $\mathcal{S}_D(\varphi) = D$. □

## 5. Optimization problems on balls in Bochner spaces

Similarly to the previous section, in this section, $L_p(S; X)$ is a uniformly convex and uniformly smooth Bochner space with dual space $L_q(S; X^*)$ with origins $\theta$ and $\theta^*$, respectively, in which $1 < p, q < \infty$ satisfying $\frac{1}{p} + \frac{1}{q} = 1$. For any $r > 0$, we define the closed, open balls and the sphere in $L_p(S; X)$ with radius $r$ and with center at $v \in L_p(S; X)$, respectively, by

$$B_p(v, r) = \left\{ f \in L_p(S; X) \colon \left( \int_S \|f(s) - v(s)\|_X^p d\mu(s) \right)^{\frac{1}{p}} \leq r \right\},$$

$$B_p^o(v, r) = \left\{ f \in L_p(S; X) : \left( \int_S \|f(s) - v(s)\|_X^p d\mu(s) \right)^{\frac{1}{p}} < r \right\},$$

$$S_p(v, r) = \left\{ f \in L_p(S; X) : \left( \int_S \|f(s) - v(s)\|_X^p d\mu(s) \right)^{\frac{1}{p}} = r \right\}.$$

Then $B_p(v, r)$ is a nonempty closed and convex subset in $L_p(S; X)$. $\{B_p^o(v, r), S_p(v, r)\}$ forms a partition of $B_p(v, r)$. In particular, if $v = \theta$, then $B_p(\theta, r)$, $B_p^o(\theta, r)$ and $S_p(\theta, r)$ are denoted by $B_p(r)$, $B_p^o(r)$ and $S_p(r)$, respectively. In this section, we consider the closed balls $B_p(r)$. All results about $B_p(r)$ can be analogously extended to $B_p(v, r)$.

**Theorem 5.1**. *For any $r > 0$, we have*

(i)   $N(B_p(r)) = B_p^o(r)$;

(ii)   $O(B_p(r)) = S_p(r)$.

(iii)   *Furthermore, for any $g \in S_p(r)$,*

    (a)   $\mathcal{S}_{B_p(r)}^{-1}(g) = \overrightarrow{\theta^*, J_p(g)} = \{t J_p(g) : 0 \leq t < \infty\}$;

    (b)   $\mathcal{S}_{B_p(r)}^{-*}(g)$ *is a curve in $L_p(S; X)$ with ending point at $\Theta$ and through $g$.*

*Proof.* Proof of (i). It is easy to see that $\Theta \in N(B_p(r))$. For arbitrary $h \in B_p^o(r)$, with

$$0 < \|h\|_{L_p(S;X)} = \left( \int_S \|h(s)\|_X^p d\mu(s) \right)^{\frac{1}{p}} < r,$$

we prove $h \in N(Bp(r))$ by the following two cases.

Case 1. $\varphi \in L_q(S; X^*)$, which satisfies $\|\varphi\|_{L_q(S;X^*)} = \left( \int_S \|\varphi(s)\|_{X^*}^q d\mu(s) \right)^{\frac{1}{q}} > 0$, and

$$\langle \varphi, h \rangle = \int_S \langle \varphi(s), h(s) \rangle d\mu(s) = 0. \tag{5.1}$$

In this case, take $u = \dfrac{r}{1 + \|\varphi\|_{L_q(S;X^*)}} J_q^*(\varphi)$. We have

$$\|u\|_{L_p(S;X)} = \frac{r}{1 + \|\varphi\|_{L_q(S;X^*)}} \left\| J_q^*(\varphi) \right\|_{L_p(S;X)}$$

$$= \frac{r}{1 + \|\varphi\|_{L_q(S;X^*)}} \|\varphi\|_{L_q(S;X^*)} < r.$$

This implies $u \in B_p^o(r)$ with $0 < \|u\|_{L_p(S;X)} < r$. It reduces

$$\langle \varphi, u \rangle = \langle \varphi, \frac{r}{1+\|\varphi\|_{L_q(S;X^*)}} J_q^*(\varphi) \rangle$$

$$= \frac{r}{1+\|\varphi\|_{L_q(S;X^*)}} \langle \varphi, J_q^*(\varphi) \rangle$$

$$= \frac{r}{1+\|\varphi\|_{L_q(S;X^*)}} \|\varphi\|_{L_q(S;X^*)}^2 > 0. \tag{5.2}$$

(5.1) and (5.2) imply that

$$\varphi \notin \mathcal{S}_{B_p(r)}^{-1}(h), \text{ for any } \varphi \in L_q(S; X^*) \text{ with } \varphi \neq \Theta^* \text{ and } \langle \varphi, h \rangle = 0. \tag{5.3}$$

Case 2. $\varphi \in L_q(S; X^*)$, which satisfies $\|\varphi\|_{L_q(S;X^*)} > 0$, and

$$\langle \varphi, h \rangle = \int_S \langle \varphi(s), h(s) \rangle d\mu(s) \neq 0. \tag{5.4}$$

By $0 < \|h\|_{L_p(S;X)} < r$, there are $s$ and $t$ with $t > 1 > s > 0$ such that $0 < \|sh\|_{L_p(S;X)} < \|h\|_{L_p(S;X)} < \|th\|_{L_p(S;X)} < r$. This implies that $th$, $sh \in B_p^o(r) \subset B_p(r)$. Then, we have

$$\max\{\langle \varphi, th \rangle, \langle \varphi, sh \rangle\} > \langle \varphi, h \rangle.$$

By (5.4), this implies

$$\varphi \notin \mathcal{S}_{B_p(r)}^{-1}(h), \text{ for any } \varphi \in L_q(S; X^*) \text{ with } \varphi \neq \Theta^* \text{ and } \langle \varphi, h \rangle \neq 0. \tag{5.5}$$

By (5.3) and (5.5), for any $h \in B_p^o(r)$, with $0 < \|h\|_{L_p(S;X)} < r$, we have

$$\varphi \notin \mathcal{S}_{B_p(r)}^{-1}(h), \text{ for any } \varphi \in L_q(S; X^*) \text{ with } \varphi \neq \Theta^*.$$

This implies

$$h \in N(B_p(r)), \text{ for any } h \in B_p^o(r), \text{ with } \|h\|_{L_p(S;X)} > 0.$$

It follows that

$$B_p^o(r) \subseteq N(B_p(r)). \tag{5.6}$$

Next, we prove (ii). Then, we will complete the proof of (i).

Proof of (ii). For any $g \in S_p(r)$, we have

$$\langle J_p g, g \rangle = \|J_p g\|_{L_q(S;X^*)}^2 = \|g\|_{L_p(S;X)}^2 = r^2.$$

This implies

$$\langle J_p g, f \rangle \leq \|J_p g\|_{L_q(S;X^*)} \|f\|_{L_p(S;X)} \leq r^2 = \langle J_p g, g \rangle, \text{ for every } f \in B_p(r). \tag{5.7}$$

Then we have

$$J_p g \in \mathcal{S}_{B_p(r)}^{-1}(g), \text{ for any } g \in S_p(r).$$

By (5.7), this implies

$$\{tJ_p(g): 0 \leq t < \infty\} \subseteq \mathcal{S}_{B_p(r)}^{-1}(g), \text{ for any } g \in S_p(r). \qquad (5.8)$$

On the other hand, for any $\psi \in \mathcal{S}_{B_p(r)}^{-1}(g)$ with $\psi \neq \theta^*$, we have

$$\langle \psi, f \rangle \leq \langle \psi, g \rangle, \text{ for every } f \in B_p(r).$$

For $0 \leq t < \infty$, it must satisfy

$$\langle t\psi, f \rangle \leq \langle t\psi, g \rangle, \text{ for every } f \in B_p(r).$$

So, we can assume that $\|\psi\|_{L_q(S;X^*)} = r$. This implies

$$\left\| J_q^* \psi \right\|_{L_p(S;X)} = \|\psi\|_{L_q(S;X^*)} = r.$$

By $\psi \in \mathcal{S}_{B_p(r)}^{-1}(g)$ and by $J_q^* \psi \in B_p(r)$, we have

$$r^2 = \|\psi\|_{L_q(S;X^*)} \|g\|_{L_p(S;X)} \geq \langle \psi, g \rangle \geq \langle \psi, J_q^* \psi \rangle = r^2.$$

This implies $\langle \psi, g \rangle = \|\psi\|_{L_q(S;X^*)} \|g\|_{L_p(S;X)}$. Since $\|\psi\|_{L_q(S;X^*)} = \|g\|_{L_p(S;X)} = r$, we obtain

$$\psi = J_p(g), \text{ for any } \psi \in \mathcal{S}_{B_p(r)}^{-1}(g) \text{ with } \|\psi\|_{L_q(S;X^*)} = r.$$

This implies

$$\mathcal{S}_{B_p(r)}^{-1}(g) \subseteq \{tJ_p(g): 0 \leq t < \infty\}. \qquad (5.9)$$

By (5.8) and (5.9), we obtain

$$\mathcal{S}_{B_p(r)}^{-1}(g) = \{tJ_p(g): 0 \leq t < \infty\}.$$

This proves (a) of part (ii). It follows

$$S_p(r) \subseteq O(B_p(r)). \qquad (5.10)$$

Since $\{N(B_p(r)), O(B_p(r))\}$ forms a partition of $B_p(r)$, then by (5.6) and (5.10), parts (i) and (ii) are simultaneously proved. Since

$$\mathcal{S}_{B_p(r)}^{-*}(g) = J_q^*(\mathcal{S}_{B_p(r)}^{-1}(g)),$$

the mapping $J_q^*$ is a one to one, onto and continuous mapping satisfying $J_q^* J_p = I_{L_p(S;X)}$, then, part (b) of (ii) follows from (a) immediately. □

**Corollary 5.2.** *Let $r > 0$ be given. Let $g = \sum_{i=1}^n (1_{A_i} \otimes x_i)$ be a $\mu$-simple function in $L_p(S; X)$ with $\|g\|_{L_p(S;X)} = r$. then*

$$\mathcal{S}_{B_p(r)}^{-1}(g) = \{t\sum_{i=1}^n \|x_i\|_X^{p-2}(1_{A_i} \otimes x_i^*): 0 \leq t < \infty\}. \qquad (5.11)$$

*Here, $x_i^* = J_X x_i$, for $i = 1, 2, \ldots, n$.*

*Proof.* By Corollary 3.3 in [13], we have

$$J_p\left(\sum_{i=1}^n (1_{A_i} \otimes x_i)\right)$$

$$= \frac{1}{\left(\sum_{j=1}^n \|x_j\|_X^p \mu(A_j)\right)^{\frac{1}{q} - \frac{1}{p}}} \sum_{i=1}^n \|x_i\|_X^{p-2} (1_{A_i} \otimes J_X x_i)$$

$$= \frac{1}{r^{p\left(\frac{1}{q} - \frac{1}{p}\right)}} \sum_{i=1}^n \|x_i\|_X^{p-2} (1_{A_i} \otimes J_X x_i).$$

By (a) of part (ii) in Theorem 5.1, this implies (5.11). □

In particular, when $n = 1$ in Corollary 5.2, we obtain

**Corollary 5.3.** *Let $r > 0$ be given. Let $A \in \mathcal{A}$ with $0 < \mu(A) < \infty$, $x \in X$ with $\mu(A)^{\frac{1}{p}} \|x\|_X = r$. Then, we have*

$$\mathcal{S}_{B_p(r)}^{-1}(1_A \otimes x) = \{t(1_A \otimes x^*): 0 \le t < \infty\}.$$

*Here, $x^* = J_X x$.*

*Proof.* Since $\|1_A \otimes x\|_{L_p(S; X)} = \mu(A)^{\frac{1}{p}} \|x\|_X = r$. This corollary follows from Corollary 5.3 immediately, in which $n = 1$. □

## 6. Analytic representations of projections on balls in Bochner spaces

In the theories of optimization, approximation, variational inequality and fixed point in Banach spaces, the three projections $P_C$ (the standard metric projection), $\pi_C$ (the generalized projection) and $\Pi_C$ (the generalized metric projection) play very important roles for the proofs of existence theorem, solution approximations in these fields. In particular, they also play crucial roles in optimization theory in Bochner spaces. In section 2, we recalled some properties of the three projections. For example, we have the connections between these three projections: $\Pi_C \ne P_C$, in general; and when $X$ is a Hilbert space, then $\pi_C = \Pi_C = P_C$. It is evident to see that the connections of these projections depend on the considered subset $C$. In this section, we consider the case that $C$ is a closed ball in some Bochner spaces and we will present the analytic representations for the projections $\pi_C$, $\Pi_C$ and $P_C$.

Similarly to the previous sections, in this section, let $L_p(S; X)$ be a uniformly convex and uniformly smooth Bochner space with dual space $L_q(S; X^*)$ with origins $\theta$ and $\theta^*$, respectively, in which $1 < p, q < \infty$ satisfying $\frac{1}{p} + \frac{1}{q} = 1$. For any $r > 0$, we denote the closed, open balls and the sphere in $L_p(S; X)$ with radius $r$ and with center at the $\theta$, by $B_p(r)$, $B_p^o(r)$ and $S_p(r)$, respectively. It is worth to note that the analytic representations for the projections $\pi_{B_p(r)}$, $\Pi_{B_p(r)}$ and $P_{B_p(r)}$ on $B_p(r)$ proved in Theorem 6.1 below can be extended to the corresponding projections on $B_p(v, r)$.

**Theorem 6.1.** *Let $L_p(S; X)$ be a uniformly convex and uniformly smooth Bochner space with dual*

space $L_q(S; X^*)$. For any $r > 0$, $g \in L_p(S; X)$ and $\varphi \in L_q(S; X^*)$, we have

(a) $\pi_{B_p(r)}(\varphi) = J_q^*(\varphi)$, if $\|\varphi\|_{L_q(S;X^*)} \leq r$;

(b) $\pi_{B_p(r)}(\varphi) = \dfrac{r}{\|\varphi\|_{L_q(S;X^*)}} J_q^*(\varphi)$, if $\|\varphi\|_{L_q(S;X^*)} > r$;

(c) $\Pi_{B_p(r)}(g) = g$, for $g \in B_p(r)$;

(d) $\Pi_{B_p(r)}(g) = \dfrac{r}{\|f\|_{L_p(S;X)}} g$, for $g \in L_p(S; X) \backslash B_p(r)$;

(e) $P_{B_p(r)}(g) = g$, for $g \in B_p(r)$;

(f) $P_{B_p(r)}(g) = \dfrac{r}{\|f\|_{L_p(S;X)}} g$, for $g \in L_p(S; X) \backslash B_p(r)$.

*Proof.* (a). As what is recalled in section 2, the generalized projection $\pi_{B_p(r)} : L_q(S; X^*) \to B_p(r)$ is defined by the Lyapunov functional $V : L_q(S; X^*) \times L_p(S; X) \to \mathbb{R}_+$:

$$V(\varphi, \pi_{B_p(r)}(\varphi)) = \inf_{f \in B_p(r)} V(\varphi, f), \text{ for any } \varphi \in L_q(S; X^*). \quad (6.1)$$

For $\varphi \in L_q(S; X^*)$ with $\|\varphi\|_{L_q(S;X^*)} \leq r$, for any $f \in B_p(r)$, we calculate

$$
\begin{aligned}
& V(\varphi, f) - V(\varphi, J_q^*(\varphi)) \\
&= (\|\varphi\|_{L_q(S;X^*)}^2 - 2\langle \varphi, f \rangle + \|f\|_{L_p(S;X)}^2) - (\|\varphi\|_{L_q(S;X^*)}^2 - 2\langle \varphi, J_q^*(\varphi) \rangle + \left\|J_q^*(\varphi)\right\|_{L_p(S;X)}^2) \\
&= -2\langle \varphi, f \rangle + \|f\|_{L_p(S;X)}^2 + \|\varphi\|_{L_q(S;X^*)}^2 \\
&\geq \|\varphi\|_{L_q(S;X^*)}^2 - 2\|\varphi\|_{L_q(S;X^*)} \|f\|_{L_p(S;X)} + \|f\|_{L_p(S;X)}^2 \\
&= (\|\varphi\|_{L_q(S;X^*)} - \|f\|_{L_p(S;X)})^2 \\
&\geq 0, \text{ for any } f \in B_p(r).
\end{aligned}
$$

Since $B_p(r)$ is a closed and convex subset of the uniformly convex and uniformly smooth Banach space $L_p(S; X)$, then $\pi_{B_p(r)} : L_q(S; X^*) \to B_p(r)$ is a one to one and onto mapping. By (6.1), the above inequalities imply that $\pi_{B_p(r)}(\varphi) = J_q^*(\varphi)$, for $\varphi \in L_q(S; X^*)$ with $\|\varphi\|_{L_q(S;X^*)} \leq r$. So, part (a) is proved.

Proof of (b). For $\varphi \in L_q(S; X^*)$ with $\|\varphi\|_{L_q(S;X^*)} > r$, for any $f \in B_p(r)$, we calculate

$$
\begin{aligned}
& V(\varphi, f) - V(\varphi, \tfrac{r J_q^*(\varphi)}{\|\varphi\|_{L_q(S;X^*)}}) \\
&= \left(\|\varphi\|_{L_q(S;X^*)}^2 - 2\langle \varphi, f \rangle + \|f\|_{L_p(S;X)}^2\right) - \left(\|\varphi\|_{L_q(S;X^*)}^2 - 2\langle \varphi, \tfrac{r J_q^*(\varphi)}{\|\varphi\|_{L_q(S;X^*)}} \rangle + \left\|\tfrac{r J_q^*(\varphi)}{\|\varphi\|_{L_q(S;X^*)}}\right\|_{L_p(S;X)}^2\right) \\
&= -2\langle \varphi, f \rangle + \|f\|_{L_p(S;X)}^2 + 2r\|\varphi\|_{L_q(S;X^*)} - r^2 \\
&\geq -2\|\varphi\|_{L_q(S;X^*)}\|f\|_{L_p(S;X)} + \|f\|_{L_p(S;X)}^2 + 2r\|\varphi\|_{L_q(S;X^*)} - r^2 \\
&= (r - \|f\|_{L_p(S;X)})\Big(2\|\varphi\|_{L_q(S;X^*)} - \|f\|_{L_p(S;X)} - r\Big)
\end{aligned}
$$

$\geq 0$, for any $f \in B_p(r)$.

This proves (b). Next, we prove (c). Notice that

$$\left\| J_p g \right\|_{L_q(S;X^*)} = \|g\|_{L_p(S;X)} \leq r, \text{ for any } g \in B_p(r). \tag{6.2}$$

Then, by (a), we have

$$\Pi_{B_p(r)}(g) = \pi_{B_p(r)}(J_p g) = J_q^*(J_p g) = g, \text{ for any } g \in B_p(r).$$

Proof of (d). Similar to (6.2), we have

$$\left\| J_p g \right\|_{L_q(S;X^*)} = \|g\|_{L_p(S;X)} > r, \text{ for any } g \in L_p(S;X) \backslash B_p(r).$$

Then by (b), we obtain

$$\Pi_{B_p(r)}(g) = \pi_{B_p(r)}(J_p g) = \frac{r}{\|J_p f\|_{L_q(S;X^*)}} J_q^*(J_p g) = \frac{r}{\|f\|_{L_p(S;X)}} g, \text{ for } g \in L_p(S;X) \backslash B_p(r).$$

Part (e) is clear. We prove part (f). For any given $g \in L_p(S;X) \backslash B_p(r)$, it satisfies $\|g\|_{L_p(S;X)} > r$. Then, for any $f \in B_p(r)$, we have

$$\begin{aligned}
&\left\langle J_p \left( g - \frac{r}{\|g\|_{L_p(S;X)}} g \right), \ \frac{r}{\|g\|_{L_p(S;X)}} g - f \right\rangle \\
&= \left( 1 - \frac{r}{\|g\|_{L_p(S;X)}} \right) \left\langle J_p(g), \ \frac{r}{\|g\|_{L_p(S;X)}} g - f \right\rangle \\
&= \left( 1 - \frac{r}{\|g\|_{L_p(S;X)}} \right) \left( r \|g\|_{L_p(S;X)} - \langle J_p(g), f \rangle \right) \\
&\geq \left( \|g\|_{L_p(S;X)} - r \right) \left( r - \|f\|_{L_p(S;X)} \right) \\
&\geq 0, \text{ for all } f \in B_p(r).
\end{aligned}$$

Since $\frac{r}{\|g\|_{L_p(S;X)}} g \in S_p(r)$, by the basic variational principle of $P_{B_p(r)}$, this implies

$$P_{B_p(r)}(g) = \frac{r}{\|g\|_{L_p(S;X)}} g, \text{ for } g \in L_p(S;X) \backslash B_p(r). \qquad \square$$

If we apply Theorem 6.1 to simple functions, we have the following results.

**Corollary 6.2**. *Let* $A \in \mathcal{A}$ *with* $\mu(A) = 1$, $x \in X$ *and* $y^* \in X^*$. *For any* $r > 0$, *we have*

(a) $\pi_{B_p(r)}(1_A \otimes y^*) = 1_A \otimes J_{X^*}(y^*), \ \text{if } \|y^*\|_{X^*} \leq r$;

(b) $\pi_{B_p(r)}(1_A \otimes y^*) = \frac{r}{\|y^*\|_{X^*}} (1_A \otimes J_{X^*}(y^*)), \ \text{if } \|y^*\|_{X^*} > r$;

(c) $\Pi_{B_p(r)}(1_A \otimes x) = 1_A \otimes x, \ \text{if } \|x\|_X \leq r$;

(d) $\Pi_{B_p(r)}(1_A \otimes x) = \frac{r}{\|x\|_X} (1_A \otimes x), \ \text{if } \|x\|_X > r$;

(e) $P_{B_p(r)}(1_A \otimes x) = 1_A \otimes x, \ \text{if } \|x\|_X \leq r$;

(f) $P_{B_p(r)}(1_A \otimes x) = \frac{r}{\|x\|_X}(1_A \otimes x),\ if\ \|x\|_X > r.$

*Proof.* By $\mu(A) = 1$, we calculate

$$\|1_A \otimes y^*\|_{L_q(S;X^*)} = \|y^*\|_{X^*} \quad \text{and} \quad \|1_A \otimes x\|_{L_p(S;X)} = \|x\|_X.$$

By (2.5), this corollary follows from Theorem 6.1 immediately. □

## Acknowledgments

The authors sincerely thank Professors Robert Mendris and Jen-Chih Yao for their valuable suggestions and comments, which improved this paper.

## References

[1]  Alber, Ya., Metric and generalized projection operators in Banach spaces: properties and applications, in "Theory and Applications of Nonlinear Operators of Accretive and Monotone Type" (A. Kartsatos, Ed.), Marcel Dekker, inc. (1996) 15–50.

[2]  Alber, Ya., and Yao, J. C., On the projection dynamical systems in Banach spaces, Taiwanese Journal of Mathematics, Vol. 11, No. 3 (2007) 819–847.

[3]  Barbu, Viorel and Michael R¨ockner, Stochastic variational inequalities and applications the total variation flow perturbed by linear multiplicative noise, arXiv:1209.0351v1 [math.PR] 3 Sep (2012).

[4]  Bochner, S. and Taylor, A. E., Linear functional} on certain spaces of abstractly-valued functions, Ann. of Math. 39 (1938) 913–944.

[5] Day, M. M., *Normed linear spaces*, 2nd rev. ed., Academic Press, New York; Springer-Verlag, Berlin, MR 26 #2847 (1962).

[6]  Dowling, P. N., Hu, Z., and Mupasiri, D., Complex convexity in Lebesgue-Bochner function spaces Trans. Amer. Math. Soc. Volume 348 Number 1 (1996).

[7]  Dunford N. and Schwartz, J. T., Linear operators. I: General theory, Pure and Appl. Math., vol. 7, Interscience, New York, (1958).

[8]  Elsanousi, S. A., Random variational inequalities in the absence of sample path continuity, Journal of Applied Mathematics and Stochastic Analysis 8, Number 1 (1995) 91-97.

[9]  Griewank Andreas Graduiertenkolleg, *Function space optimization*, 1128 Institut f`ur Mathematik Humboldt-Universit¨at zu Berlin Script, version June 4, 2008 By Levis Eneya & Lutz Lehmann Winter Semester 2007/2008

[10] Hytönen, T., Neerven, J., Veraar, M. and Weis, L., *Analysis in Banach spaces*, Ergebnisse der Mathematik und ihrer Grenzgebiete. 3. Folge / A Series of Modern Surveys in Mathematics, Volume 63, Springer International Publishing AG (2016).

[11] Ito, Kazufumi, *Functional Analysis and Optimization*, Department of Mathematics, North Carolina State University, Raleigh, North Carolina, USA November 29 (2016).

[12] Khan, A. A., Li, J. L. and Reich, S., Generalized projection operators on general Banach spaces, to appear in *Journal of Nonlinear and Convex Analysis*.

[13] Khan, A. A., Li, J. L., and Tammer Ch., Normalized Duality Mappings in Bochner Spaces, in preparation.

[14] Khan, A. A., Li, J. L. Approximating Properties of Metric and Generalized Metric Projections in Uniformly Convex and Uniformly Smooth Banach Spaces, in preparation.

[15] Kosmol, Peter and Dieter Müller-Wichards *Optimization in Function Spaces With Stability Considerations in Orlicz Spaces*, Published by De Gruyter, Volume 13 in the series Nonlinear Analysis and Applications  https://doi.org/10.1515/9783110250213 (2011)

[16] Leonard, E. and Sundaresan, K., Geometry of Lebesgue-Bochner Function Spaces-Smoothness, Transactions of the American Mathematical Society, Vol. 198, (1974) 229–251.

[17] Li, J. L., On Geometric Properties of Bochner Spaces and Expected Fixed Point in Banach Spaces, to appear in *Journal of Nonlinear and Convex Analysis*.

[18] Li, J. L., "The Generalized Projection Operator on Reflexive Banach Spaces and Its Applications," *Journal of Mathematical Analysis and Applications*, Vol. 306 (2005) 55–71.

[19] Long, X. J., He, Y. H., and Huang N. J., Variance-based Bregman extragradient algorithm with line search for solving stochastic variational inequalities, arXiv:2208.14069v1 [math.OC]  30 Aug (2022)

[20] Mikusiński, J., *The Bochner Integral*, Mathematische Reihe (1978).

[21] Sasane, Amol, *Optimization in Function Spaces* (Aurora: Dover Modern Math Originals) First Edition (2016)

[22] Shanbhag, Uday V., Stochastic Variational Inequality Problems: Applications, Analysis, and Algorithms. INFORMS TutORials in Operations Research. Published online: 14 Oct 2014; 71-107. https://doi.org/10.1287/educ.2013.0120

[23] Sobolev, V. I., "Bochner integral", Encyclopedia of Mathematics, EMS Press (2001).

[24]  Sundaresan, K., The Radon-Nikodým Theorem for Lebesgue-Bochner function spaces, J. Funct. Anal. 24, (1977) 276–279.

[35] Takahashi, W., *Nonlinear Functional Analysis,* Yokohama Publishers, 2000.

[26] Zaslavski, Alexander J., *Optimization in Banach Spaces*, Springer Briefs in Optimization ISSN 2190-8354 ISSN 2191-575X (electronic) ISBN 978-3-031-12643-7 ISBN 978-3-031-12644-4 (eBook) https://doi.org/10.1007/978-3-031-12644-4